\documentclass[graybox]{svmult}

\usepackage{graphicx}

\usepackage{amsmath}
\usepackage{color}

\usepackage[latin1]{inputenc}
\usepackage{amsfonts}
\usepackage{tikz} 
\usetikzlibrary{matrix}
\usepackage[english]{babel}
\usepackage[latin1]{inputenc}
\usepackage[T1]{fontenc}
\usepackage{ae}

\usepackage{enumitem}

\usepackage{url}

\usepackage{algorithm}
\usepackage{algorithmicx}
\usepackage{algpseudocode}

\newenvironment{psmallmatrix}
  {\left(\begin{smallmatrix}}
  {\end{smallmatrix}\right)}

\allowdisplaybreaks[4]

\usepackage{color}
\usepackage{tikz}
\usetikzlibrary{matrix}
\usepackage{cite}           
\usepackage{mathptmx}       
\usepackage{helvet}         
\usepackage{courier}        
\usepackage{type1cm}        
\usepackage{makeidx}         
\usepackage{graphicx}        
\usepackage{multicol}        
\usepackage[bottom]{footmisc}

\newcommand{\OOtilde}{\smash{\widetilde{O}}}

\makeindex             

\begin{document}


\title{Numerical Evaluation of Elliptic Functions, Elliptic Integrals and Modular Forms}

\author{Fredrik Johansson}
\institute{Fredrik Johansson,\\
INRIA -- LFANT, CNRS -- IMB -- UMR 5251, Universit\'{e} de Bordeaux, 33400 Talence, France,
\email{fredrik.johansson@gmail.com}
}
%
%
\maketitle

\abstract{We describe algorithms to compute elliptic functions
and their relatives (Jacobi theta functions, modular forms, elliptic integrals,
and the arithmetic-geometric mean)
numerically to arbitrary precision with rigorous error bounds
for arbitrary complex variables. Implementations in
ball arithmetic are available
in the open source Arb library. We discuss the algorithms
from a concrete implementation point of view,
with focus on performance at tens to thousands of digits of precision.}

%
%

\section{Introduction}
\label{sec:1}

\vspace{1mm}\noindent
The elliptic functions and their relatives
have many applications in
mathematical physics and number theory.
Among the elliptic family of special functions,
we count the elliptic functions proper (i.e.\ doubly periodic meromorphic functions)
as well as the quasiperiodic Jacobi theta functions,
the closely related classical modular forms and modular functions
on the upper half plane,
and elliptic integrals which are the inverse functions
of elliptic functions.

\index{elliptic function}
\index{Jacobi theta function}
\index{modular form}
\index{elliptic integral}
\index{arithmetic-geometric mean}
\index{arbitrary-precision arithmetic}
\index{numerical evaluation}

Our goal is to give a modern treatment of
numerical evaluation of these functions, using
algorithms that meet several criteria:

\begin{itemize}

\item \emph{Full domain}. We should be able to compute the functions for
arbitrary complex values of all parameters where this is reasonable, with
sensible handling of branch cuts for multivalued functions.

\item \emph{Arbitrary precision}.
Precision much higher than 16-digit (or 53-bit) machine arithmetic is
sometimes needed for solving numerically ill-conditioned problems.
For example, extremely high precision evaluations are employed
in mathematical physics to find
closed-form solutions for sums and integrals using
integer relation methods~\cite{bailey2015high}. \index{integer relation}
Computations with elliptic functions and modular forms
requiring hundreds or thousands of digits are commonplace in
algebraic and analytic number theory,
for instance in the construction of discrete data such as class polynomials
from numerical approximations~\cite{enge2009complexity}. \index{class polynomial}

\item \emph{Rigorous error bounds}.
When we compute an approximation $y \approx f(x)$, we should
also be able to offer a bound $|y - f(x)| \le \varepsilon$,
accounting for all intermediate rounding and approximation errors
in the algorithm.

\item \emph{Efficiency}.
The algorithms should be robust and efficient
for arguments that are very small, large, or close to singularities.
For arbitrary-precision implementations,
a central concern is to ensure that the computational complexity
as a function of the precision does not grow too quickly.
At the same time, we must have in mind that
the algorithms with the best theoretical asymptotic
complexity are not necessarily the best in practice, and
we should not sacrifice efficiency at moderate precision (tens to hundreds of digits).

\end{itemize}

It turns out that these goals can be achieved simultaneously
and with reasonable implementation effort
thanks to the remarkable amount of structure in the
elliptic function family
--- in contrast to many other problems in numerical analysis!

The present author has implemented
routines for elliptic functions, Jacobi theta functions,
elliptic integrals
and commonly-used modular forms and functions as part of the open source Arb library for arbitrary-precision \index{Arb library} \index{ball arithmetic}
ball arithmetic~\cite{Johansson2017arb}.\footnote{Available at \url{http://arblib.org}. The functionality for modular forms and elliptic
functions can be found in the \texttt{acb\_modular} (\url{http://arblib.org/acb_modular.html})
and \texttt{acb\_elliptic} (\url{http://arblib.org/acb_elliptic.html}) modules.}
The idea behind ball arithmetic is to represent numerical approximations
with error bounds attached, as in
$\pi \in [3.14159265358979 \pm 3.57 \cdot 10^{-15}]$.
The algorithms use ball arithmetic internally for automatic
propagation of error bounds, in combination
with some pen-and-paper bounds mainly for truncations of infinite series.\footnote{Of course, for applications that do not require rigorous error bounds,
all the algorithms can just as well be implemented in ordinary floating-point
arithmetic.}

The functions in Arb can be used directly in C or via
the high-level wrappers in Sage~\cite{sagemath} or the Julia package Nemo~\cite{Fieker2017nemo}.
As an example, we use the Arb interface in Sage to evaluate the \index{Weierstrass elliptic function}
Weierstrass elliptic function $\wp$ on the lattice \index{lattice}
$(1, \tau)$ with $\tau = \tfrac{1}{2}(1 + \sqrt{3} i)$.
We check $\wp(z) = \wp(z + 5 + 6\tau)$ at the
arbitrarily chosen point $z = 2 + 2i$, here using 100-bit precision:

\begin{footnotesize}
\begin{verbatim}
sage: C = ComplexBallField(100)
sage: tau = (1 + C(-3).sqrt())/2
sage: z = C(2 + 2*I)
sage: z.elliptic_p(tau)
[-13.7772161934928750714214345 +/- 6.41e-26] + [+/- 3.51e-26]*I
sage: (z + 5 + 6*tau).elliptic_p(tau)
[-13.777216193492875071421435 +/- 9.69e-25] + [+/- 4.94e-25]*I
\end{verbatim}
\end{footnotesize}

This text covers the algorithms used in Arb
and discusses some of the implementation aspects.
The algorithms are general and work well in most situations.
However, we note that the code in Arb does not use the best available algorithms
in all cases, and we will point out some
of the possible improvements.

There is a vast body of literature on elliptic functions and integrals,
and we will not be able to explore the full breadth of
computational approaches.
We will, in particular, focus on arbitrary-precision arithmetic
and omit techniques that only matter in machine precision.
A good overview and a comprehensive bibliography can be found in chapters 19, 20, 22 and 23 of the
NIST Handbook of Mathematical Functions~\cite{Olver2010}
or its online counterpart, the Digital Library of Mathematical Functions\footnote{\url{https://dlmf.nist.gov/}}.
Cohen's book on computational number theory \cite{cohen2013course} is also a useful resource.

Many other packages and computer algebra systems
also provide good support for evaluating elliptic and related functions,
though not with rigorous error bounds; we mention Pari/GP~\cite{PARI2} \index{Pari/GP}
and of course Maple and Mathematica.
For a nice application of Weierstrass elliptic functions in astrodynamics \index{Weierstrass elliptic function}
and a fast machine-precision implementation of these functions, we mention
the work by Izzo and Biscani~\cite{izzo2016astrodynamics}.

The algorithms that we review are well known, but
they are sometimes described without discussing arbitrary complex variables,
variable precision, or error bounds.
We attempt to provide an account that is complementary to the
existing literature, and we also
discuss some minor improvements to the algorithms as they
are usually presented. For example, we have optimized Carlson's algorithm \index{Carlson symmetric form}
for symmetric elliptic integrals to reduce the asymptotic complexity
at high precision (section \ref{sect:rfalg}),
and we make several observations
about the deployment of ball arithmetic.

\section{General strategy}

\vspace{1mm}\noindent
Algorithms for evaluating mathematical functions often
have two stages: argument reduction,
followed by evaluation of a series expansion \cite{brent2010modern,muller2006elementary}.
\index{argument reduction}

Minimax polynomial or rational function
approximations are usually preferred for univariate functions
in machine precision, but truncated
Taylor series expansions are the tools of choice in arbitrary-precision arithmetic,
for two reasons. First, precomputing minimax approximations is not practical,
and second, we can exploit the
fact that the polynomials arising from series expansions of special
functions are typically not of generic type but highly structured.

\index{numerical stability}

Argument reduction consists of applying functional equations
to move the argument to a part of the domain where the series
expansion converges faster.
In many cases, argument reduction is needed to ensure convergence in the first place.
Argument reduction also tends to improve numerical stability,
in particular by avoiding alternating series with large terms
that would lead to catastrophic cancellation.

\begin{table}[t]
\begin{center}
\caption{\label{tab:algoverview}Methods for computation of elliptic functions and integrals.
This table illustrates the analogies between the elliptic function and elliptic integral cases,
and the simplifications between the general (arbitrary parameters) and special (some parameters fixed at special values) cases.}
\begin{small}
\renewcommand{\arraystretch}{1.3}
\setlength{\tabcolsep}{.5em}
\begin{tabular}{|p{2.4cm} | p{3.8cm} | p{4.3cm}|} \hline
              & \textbf{Elliptic functions} & \textbf{Elliptic integrals} \\ \hline \hline
\textit{General case}  & \textit{\mbox{Elliptic functions}, \mbox{Jacobi theta functions}} & \textit{Incomplete elliptic integrals} \\ \hline
Argument reduction  & \mbox{Reduction to standard domain} (\mbox{modular transformations}, \mbox{periodicity)} & \mbox{Contraction of parameters}
    (\mbox{linear symmetric} \mbox{transformations}) \\
Series expansions   & Theta function $q$-series & \mbox{Multivariate hypergeometric} \mbox{series} \\ \hline \hline
\textit{Special case}  & \textit{\mbox{Modular forms \& functions,} \mbox{theta constants}} & \textit{\mbox{Complete elliptic integrals}, arithmetic-geometric mean} \\ \hline
Argument reduction  & \mbox{Reduction to standard domain} (modular transformations)  & \mbox{Contraction of parameters} (quadratic transformations) \\
Series expansions   & \mbox{Theta constant and eta function} \mbox{$q$-series} & \mbox{Classical ${}_2F_1$ hypergeometric} \mbox{series} \\ \hline
\end{tabular}
\end{small}
\end{center}
\end{table}

\index{elliptic integral}
\index{complete elliptic integral}
\index{elliptic function}
\index{modular form}
\index{modular function}
\index{incomplete elliptic integral}
\index{hypergeometric series}
\index{Gauss hypergeometric function}

The classical elliptic and modular functions are no exception to this
general pattern, as shown in Table~\ref{tab:algoverview}.
For elliptic integrals, the argument reduction consists of using
contracting transformations to reduce the distance
between the function arguments,
and the series expansions are hypergeometric series (in one
or several variables).
For the elliptic and modular functions, the argument reduction
consists of using modular transformations and periodicity to move the
lattice parameter to the fundamental domain and the argument to a lattice cell near the origin, and the
series expansions are the sparse $q$-series of Jacobi theta functions.

\index{modular transformation}
\index{$q$-series}
\index{Jacobi theta function}
\index{hypergeometric series}
\index{modular form}
\index{modular function}

In the following text, we will first discuss the computation
of elliptic functions starting with the special case of modular forms
and functions before turning to general elliptic and Jacobi theta functions.
Then, we discuss elliptic integrals,
first covering the easier case of complete integrals
before concluding with the treatment of incomplete integrals.

We comment briefly on error bounds. Since ball arithmetic automatically
tracks the error propagation 
during series evaluation and through argument reduction steps,
the only error analysis that needs to be done by hand
is to bound the series truncation errors.
If $f(x) = \sum_{k=0}^{\infty} t_k(x)$,
we compute $\sum_{k=0}^N t_k(x)$ and then add 
the ball $[\pm \varepsilon]$
or $[\pm \varepsilon] + [\pm \varepsilon]i$
where $\varepsilon$ is an upper bound for $|R_N(x)| = |\sum_{k=N+1}^{\infty} t_k(x)|$.
Such a bound is often readily obtained by comparison with a geometric series,
i.e.\ if $|t_k(x)| \le A C^k$ with $0 \le C < 1$, then $|R_N(x)| \le \sum_{k=N+1} A C^k = A C^N / (1-C)$.
In some cases, further error analysis can be useful to improve
the quality (tightness) of the ball enclosures.

\index{asymptotic complexity}

For arbitrary-precision evaluation, we wish to minimize
the computational complexity as a function of the precision $p$.
The complexity is often measured by counting arithmetic operations.
The actual time complexity must account for the fact that
arithmetic operations have a bit complexity of $\OOtilde(p)$ (where the
$\OOtilde$ notation ignores logarithmic factors).
In some situations, it is better to use a model of complexity
that distinguishes between ``scalar'' arithmetic operations
(such as addition of
two $p$-bit numbers or multiplication of a $p$-bit number by a small integer)
and ``nonscalar'' arithmetic operations (such as multiplication
of two general $p$-bit numbers).

\subsection{The exponential function}

\vspace{1mm}\noindent
We illustrate these principles with a commonly used algorithm to compute
the exponential function $e^x$ of a real argument $x$
to $p$-bit precision.

\begin{itemize}[listparindent=1.5em]

\index{exponential function}
\index{argument reduction}

\item \emph{Argument reduction.}
We first use $e^x = 2^n e^t$ with
$t = x - n \log(2)$ and $n = \lfloor x / \log(2) \rfloor$ which
ensures that $t \in [0,\log(2))$.
At this point, the usual Taylor series $e^t = 1 + t + \tfrac{1}{2} t^2 + \ldots$
does not suffer from cancellation, and we only need $O(p/\log p)$
terms for a relative error of $2^{-p}$ independent of the initial size of~$|x|$.
As a second argument reduction step, we write
$e^t = (e^u)^{2^r}$ with $u = t/2^r$, which reduces the number $N$
of needed Taylor series terms to $O(p/r)$.

Balancing $N = O(p/r)$
against the number $r$ of squarings needed to reconstruct
$e^t$ from $e^u$, it is optimal to choose $r \approx p^{0.5}$. This
gives an algorithm for $e^x$ requiring $O(p^{0.5})$ arithmetic operations
on $p$-bit numbers,
which translates to a time complexity of $\OOtilde(p^{1.5})$.

\item \emph{Series evaluation.}
As an additional improvement, we can exploit the structure of
the Taylor series of the exponential function.
For example, $\sum_{k=0}^8 \tfrac{1}{k!} x^k$ can be evaluated as
\begin{equation}
1 + x + \tfrac{1}{2} \left(
x^2 + \tfrac{1}{3} x^3 \left(
1 + \tfrac{1}{4} \left(
x + \tfrac{1}{5} \left(
x^2 + \tfrac{1}{6} x^3 \left(
1 + \tfrac{1}{7} \left(
x + \tfrac{1}{8} x^2
\right) \right) \right) \right) \right) \right)
\label{eq:rectsplit}
\end{equation}
where we have extracted the power $x^3$ repeatedly and used the
fact that the ratios between successive
coefficients are small integers.
As a result, we only need four nonscalar multiplications
involving $x$ (to compute $x^2, x^3$, and for the
two multiplications by $x^3$), while the remaining
operations are scalar divisions. With further rewriting, the scalar divisions
can be replaced by even cheaper scalar multiplications.

In general, to evaluate a polynomial of degree $N$ with scalar coefficients
at a nonscalar argument $x$,
we can compute $x^2, \ldots, x^m$ once
and then use Horner's rule with respect to $x^m$, for $m \approx N^{0.5}$,
which reduces the total number of nonscalar multiplications to
about $2 N^{0.5}$~\cite{PatersonStockmeyer1973}. This trick is
sometimes called \emph{rectangular splitting}.
To motivate this terminology, picture the terms of
the polynomial laid out as a matrix with $m$ columns and $N/m$ rows.

\index{rectangular splitting}
\index{arithmetic-geometric mean}

In view of this improvement to the series evaluation, it turns out to be more efficient
in practice to choose the tuning parameter $r$ used for argument reduction
$e^t = (e^{t/2^r})^{2^r}$ slightly
smaller, say about $r \approx p^{0.4}$ for realistic $p$.
The algorithm combining optimal argument reduction with rectangular splitting
for evaluation of elementary functions such as $e^x$
is due to Smith~\cite{Smith1989}.

There are asymptotically faster algorithms that permit
evaluating elementary functions using only $O(\log p)$
arithmetic operations
(that is, in $\OOtilde(p)$ time), for instance
based on the AGM (discussed in section~\ref{sect:complagm}
for computing elliptic integrals), but Smith's algorithm
is more efficient in practice for moderate $p$,
and in some situations still wins for $p$ as large as $10^5$.

\end{itemize}

\section{Modular forms and functions}
\vspace{1mm}\noindent
A modular transformation $g$ is a linear fractional transformation
on the upper half plane $\mathbb{H} = \{\tau : \mathbb{C} : \operatorname{Im}(\tau) > 0\}$
of the form $g(\tau) = (a \tau + b) / (c \tau + d)$
where $a,b,c,d$ are integers with $ad - bc = 1$.
We can assume that $c \ge 0$.
The group of modular transformations (known as the modular group) can
be identified
with the projective special linear group
$\text{PSL}(2,\mathbb{Z})$, where $g$ is represented by 
the matrix $\begin{psmallmatrix} a & b \\ c & d \end{psmallmatrix}$
and composition corresponds to matrix multiplication.

\index{modular transformation}
\index{modular group}
\index{modular form}
\index{modular function}

A modular form (of weight $k$) is a holomorphic function on $\mathbb{H}$
satisfying the functional equation $f(g(\tau)) = (c\tau+d)^k f(\tau)$
for every modular transformation $g$,
with the additional technical requirement of being
holomorphic as $\tau \to i \infty$.
A modular form of weight $k = 0$ must be a constant function, but nontrival
solutions of the above functional equation are possible if we allow poles.
A meromorphic function on $\mathbb{H}$ satisfying $f(g(\tau)) = f(\tau)$
is called a modular function.%

\index{Fourier series}
\index{$q$-series}

Every modular form or function is periodic with $f(\tau+1) = f(\tau)$
and has a Fourier series (or $q$-series)
\begin{equation}
f(\tau) = \sum_{n=-m}^{\infty} a_n q^n, \quad q^{2\pi i \tau}
\label{eq:qseries}
\end{equation}
where $m = 0$ in the case of a modular form.
The fundamental tool in numerical evaluation of modular forms
and functions is to evaluate a truncation of such a $q$-series.
Since $\tau$ has positive imaginary part, the quantity $q$
always satisfies $|q| < 1$, and
provided that an explicit bound for the coefficient sequence $a_n$
is known, tails of \eqref{eq:qseries} are easily bounded by a geometric series.

\subsection{Argument reduction}

\index{argument reduction}
\index{fundamental domain}

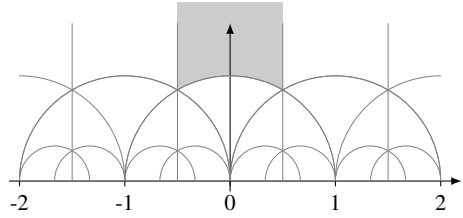
\begin{figure}
\begin{center}

\pgfdeclarelayer{background}
\pgfsetlayers{background,main}

\pgfmathsetmacro{\myxlow}{-2}
\pgfmathsetmacro{\myxhigh}{2}
\pgfmathsetmacro{\myiterations}{2}

\begin{tikzpicture}[scale=1.4]
    \draw[-latex](\myxlow-0.1,0) -- (\myxhigh+0.2,0);
    \pgfmathsetmacro{\succofmyxlow}{\myxlow+1.0}
    \draw (-2,0) -- (-2,-0.05) node[below,font=\small] {-2};
    \draw (-1,0) -- (-1,-0.05) node[below,font=\small] {-1};
    \draw (0,0) -- (0,-0.05) node[below,font=\small] {0};
    \draw (+1,0) -- (+1,-0.05) node[below,font=\small] {1};
    \draw (+2,0) -- (+2,-0.05) node[below,font=\small] {2};

    \draw[-latex](0,-0.1) -- (0,1.5);
    \draw[very thin, gray](-1.5,0) -- (-1.5,1.5);
    \draw[very thin, gray](-0.5,0) -- (-0.5,1.5);
    \draw[very thin, gray](0.5,0) -- (0.5,1.5);
    \draw[very thin, gray](1.5,0) -- (1.5,1.5);
    \begin{scope}   
        \clip (\myxlow,0) rectangle (\myxhigh,1.1);
        \foreach \i in {1,...,\myiterations}
        {   \pgfmathsetmacro{\mysecondelement}{\myxlow+1/pow(2,floor(\i/3))}
            \pgfmathsetmacro{\myradius}{pow(1/3,\i-1}
            \foreach \x in {-2,\mysecondelement,...,2}
            {   \draw[very thin, gray] (\x,0) arc(0:180:\myradius);
                \draw[very thin, gray] (\x,0) arc(180:0:\myradius);
            }   
        }
    \end{scope}
    \begin{scope}
        \begin{pgfonlayer}{background}
            \clip (-0.5,0) rectangle (0.5,1.7);
            \clip (1,1.7) -| (-1,0) arc (180:0:1) -- cycle;
            \fill[gray,opacity=0.4] (-1,-1) rectangle (1,2);
        \end{pgfonlayer}
    \end{scope}
\end{tikzpicture}
\caption{The shaded region shows the canonical fundamental domain for the action of the modular group
on the upper half plane.}
\label{fig:fundom}
\end{center}
\end{figure}

The $q$-series \eqref{eq:qseries} always converges,
but the convergence is slow for $\tau$ close to the real line
where $|q| \approx 1$.
However, we can always find a modular transformation $g$ that moves
$\tau$ to the fundamental domain
$\{\tau \in \mathbb{H} : |\tau| \ge 1, |\operatorname{Re}(\tau)| \le \tfrac{1}{2}\}$ (see Fig.~\ref{fig:fundom}).
This ensures $\smash{|q| \le e^{-\pi \sqrt{3}} \approx 0.00433}$ which
makes the convergence extremely rapid.

Technically, the fundamental domain does not include half of the boundary
(meaning that $\mathbb{H}$ is tiled by copies of the fundamental domain under
the action of the modular group), but this does not matter for the algorithm.
In fact, it is sufficient to put $\tau$ within some small
distance $\varepsilon$ of the fundamental domain, and
this relaxation is especially useful in ball arithmetic since a ball may overlap with
the boundary.

\index{ball arithmetic}

The well-known algorithm to construct $g$ (see Cohen \cite[Algorithm 7.4.2]{cohen2013course})
repeatedly applies the generators
$\tau \to \tau + 1$ and $\tau \to -1/\tau$ of the modular group:

\begin{enumerate}
\item Set $g \gets \begin{psmallmatrix} 1 & 0 \\ 0 & 1 \end{psmallmatrix}$.
\item Set $\tau \gets \tau + n$, $g \gets \begin{psmallmatrix} 1 & n \\ 0 & 1 \end{psmallmatrix} g$
where $n = -\left \lfloor \operatorname{Re}(\tau) + \tfrac{1}{2} \right \rfloor$.
\item If $|\tau| < 1 - \varepsilon$, set $\tau \gets -1/\tau$ and $g \gets \begin{psmallmatrix} 0 & -1 \\ 1 & 0 \end{psmallmatrix} g$ and go to step 2;
otherwise, stop.
\end{enumerate}

Exact integer operations should be used for the matrices,
but we can perform all operations involving $\tau$
in this algorithm
using heuristic floating-point approximations,
starting with the midpoint of the initial ball.
Once $g$ has been constructed,
we evaluate $g(\tau)$ and $f(g(\tau))$ as usual in ball arithmetic.

\index{ball arithmetic}
\index{transformation matrix}
\index{modular form}
\index{arbitrary-precision arithmetic}

Indeed, it is important to construct the transformation matrix $g$ separately and then
apply the functional equation for the modular form in a single
step rather than applying the generating transformations
iteratively in ball arithmetic.
This both serves to minimize numerical instability
and to optimize performance.
The precision needed to construct $g$ only depends on
the size of the entries of $g$ and not on the precision for evaluating $f(\tau)$.
If $\tau \approx 2^{-p} i$, we only need about $2p$ bits to construct $g$
even if the overall precision is thousands of digits.
The implementation in Arb uses virtually costless
machine floating-point arithmetic to construct $g$
when 53-bit arithmetic is sufficient, switching to arbitrary-precision arithmetic
only when necessary.

\subsection{Standard functions}

\index{Dedekind eta function}
\index{theta constant}

Several commonly-used modular forms and functions are implemented in Arb.
The basic building blocks are the Dedekind eta function
\begin{equation}
\eta(\tau) = e^{\pi i \tau/12} \sum_{n=-\infty}^{\infty} (-1)^n q^{(3n^2-n)/2}, \quad q = e^{2 \pi i \tau},
\label{eq:etaseries}
\end{equation}
and the theta constants $\theta_j \equiv \theta_j(\tau)$,
\begin{equation}
\theta_2(\tau) = e^{\pi i \tau / 4} \sum_{n=-\infty}^{\infty} q^{n(n+1)}, \quad
\theta_3(\tau) = \!\!\sum_{n=-\infty}^{\infty} q^{n^2}, \quad
\theta_4(\tau) = \!\!\sum_{n=-\infty}^{\infty} (-1)^n q^{n^2}
\label{eq:thetaseries}
\end{equation}
in which (as a potential source of confusion) $q = e^{\pi i \tau}$.

It is useful to represent other modular forms in terms of these particular functions
since their $q$-series are extremely sparse (requiring only $O(p^{0.5})$
terms for $p$-bit accuracy, which leads to $\OOtilde(p^{1.5})$ bit complexity) and only have coefficients $\pm 1$.
We give a few examples of derived functions:
\begin{itemize}
\item Modular functions are precisely the rational functions of the $j$-invariant
\begin{equation}
j(\tau) = 32 \frac{(\theta_2^8+\theta_3^8+\theta_4^8)^3}{(\theta_2 \theta_3 \theta_4)^8}, \quad j\left(\frac{a\tau+b}{c\tau+d}\right) = j(\tau).
\label{eq:modj}
\end{equation}
\item The modular discriminant is a modular form of weight 12, given by
\begin{equation}
\Delta(\tau) = \eta(\tau)^{24}, \quad \Delta\left(\frac{a\tau+b}{c\tau+d}\right) = (c\tau+d)^{12} \Delta(\tau).
\label{eq:moddelta}
\end{equation}
\item Eisenstein series are modular forms of weight $2k$ for $k \ge 2$, given by
\begin{equation}
G_{2k}(\tau) = \sum_{m^2 + n^2 \ne 0} \frac{1}{(m+n\tau )^{2k}}, \quad G_{2k} \left(\frac{a\tau+b}{c\tau+d}\right) = (c\tau+d)^{2k} G_{2k}(\tau)
\label{eq:modeisen}
\end{equation}
where we compute $G_4(\tau)$, $G_6(\tau)$ via theta constants using
$$G_4(\tau) = \frac{\pi^4}{90} \left(\theta_2^8 + \theta_3^8 + \theta_4^8\right), \quad
G_6(\tau) = \frac{\pi^6}{945} \left(-3\theta_2^8 (\theta_3^4 + \theta_4^4) + \theta_3^{12} + \theta_4^{12}\right)$$
and obtain the higher-index values using recurrence relations.
\end{itemize}

\index{modular form}
\index{modular function}
\index{$j$-invariant}
\index{modular discriminant}
\index{Eisenstein series}
\index{modular transformation}
\index{Kronecker symbol}
\index{Jacobi theta function}

The Dedekind eta function itself transforms as
\begin{equation}
\eta\left(\frac{a\tau+b}{c\tau+d}\right) = \varepsilon (a,b,c,d)
    \, \sqrt{c\tau+d} \, \eta(\tau)
\label{eq:etatransform}
\end{equation}
where $\varepsilon(a,b,c,d) = \exp(\pi i R / 12)$ is a 24th root of unity.
The integer $R$ (mod 24) can be computed using Kronecker symbols \cite[section 74]{Rademacher1973}.
The modular transformations for theta constants
are a special case of the formulas for theta functions
given below in section~\ref{sect:thetatransform}. However, we avoid
using these transformations directly when computing the functions
\eqref{eq:moddelta}, \eqref{eq:modeisen} and \eqref{eq:modj}:
it is better to apply the simpler argument reductions
for the top-level functions and then evaluate the series expansions
\eqref{eq:etaseries} or \eqref{eq:thetaseries} when $\tau$ is already reduced.

\subsection{Fast evaluation of $q$-series}

\index{$q$-series}
\index{addition sequence}
\index{rectangular splitting}

The powers of $q$ appearing in \eqref{eq:etaseries} and \eqref{eq:thetaseries}
are easily generated using two multiplications per term since
the exponents are successive values of quadratic polynomials.
The cost can nearly be halved
using short addition sequences~\cite[Algorithm 2]{enge2018short}.
The cost can be reduced even further
by combining addition sequences with rectangular splitting\cite[section 5]{enge2018short}.
Here, the idea is to factor out some power $q^m$ as in \eqref{eq:rectsplit},
but~$m$ must be chosen in a particular way --- for example, in the case of the theta series $\smash{\sum_{n=1}^{\infty} q^{n^2}}$,
$m$ is chosen so that there are few distinct quadratic residues modulo~$m$.
In Arb, these optimizations save roughly a factor four over the naive algorithm.

\section{Elliptic and theta functions}
\vspace{1mm}\noindent
An elliptic function with respect to a lattice in $\mathbb{C}$ with periods $\omega_1, \omega_2$
is a meromorphic function satisfying $f(z + m \omega_1 + n \omega_2) = f(z)$ for all $z \in \mathbb{C}$
and all $m, n \in \mathbb{Z}$.
By making a linear change of variables,
we can assume that $\omega_1 = 1$ and $\omega_2 = \tau \in \mathbb{H}$.
The elliptic functions with a fixed lattice parameter $\tau$ form a field, which is generated
by the Weierstrass elliptic function
\begin{equation}
\wp(z, \tau) = \frac{1}{z^2} + \sum_{n^2+m^2 \ne 0}
    \left[ \frac{1}{(z+m+n\tau)^2} - \frac{1}{(m+n\tau)^2} \right]
\end{equation}
together with its $z$-derivative $\wp'(z,\tau)$.

\index{elliptic function}
\index{Weierstrass elliptic function}
\index{lattice}
\index{period}
\index{lattice parameter}
\index{Jacobi theta function}

The building blocks for elliptic functions are the Jacobi theta functions
\begin{align}
\theta_1(z,\tau) &= \!\!\!\sum_{n=-\infty}^{\infty} e^{\pi i [(n + \tfrac12)^2 \tau + (2n + 1) z + n - \tfrac12]} \nonumber
                 = 2 q_{4} \sum_{n=0}^{\infty} (-1)^n q^{n(n+1)} \sin((2n+1) \pi z) \\ \nonumber
                 &= -i q_4 \sum_{n=0}^{\infty} (-1)^n q^{n(n+1)} (w^{2n+1} - v^{2n+1}), \\ \nonumber
\theta_2(z,\tau) &= \!\!\!\sum_{n=-\infty}^{\infty} e^{\pi i [(n + \tfrac12)^2 \tau + (2n + 1) z]} \nonumber
                 = 2 q_{4} \sum_{n=0}^{\infty} q^{n(n+1)} \cos((2n+1) \pi z) \\
\label{eq:genthetaseries}
                 &= q_4 \sum_{n=0}^{\infty} q^{n(n+1)} (w^{2n+1} + v^{2n+1}), \\ \nonumber
\theta_3(z,\tau) &= \!\!\!\sum_{n=-\infty}^{\infty} e^{\pi i [n^2 \tau + 2n z]}
                 = 1 + 2 \sum_{n=1}^{\infty} q^{n^2} \cos(2n \pi z) 
                 = 1 + \sum_{n=1}^{\infty} q^{n^2} (w^{2n} + v^{2n}), \\ \nonumber
\theta_4(z,\tau) &= \!\!\!\sum_{n=-\infty}^{\infty} e^{\pi i [n^2 \tau + 2n z + n]} \nonumber
                 = 1 + 2 \sum_{n=1}^{\infty} (-1)^n q^{n^2} \cos(2n \pi z) \\ \nonumber
                 &= 1 + \sum_{n=1}^{\infty} (-1)^n q^{n^2} (w^{2n} + v^{2n}), \nonumber
\end{align}
where $q = e^{\pi i \tau}$, $q_{4} = e^{\pi i \tau / 4}$, $w = e^{\pi i z}$, $v = w^{-1}$.
The theta functions
are quasielliptic functions of $z$, having period or half-period 1 and quasiperiod~$\tau$
(a shift by $\tau$ introduces an exponential prefactor).
With $z = 0$, the theta functions $\theta_2, \theta_3, \theta_4$
reduce to the corresponding theta constants,
while $\theta_1(0,\tau) = 0$ identically.

\index{quasielliptic function}
\index{Arb library}
\index{power series arithmetic}

Arb provides a complete implementation of the Jacobi theta functions themselves
as well as the Weierstrass elliptic function which is computed as
\begin{equation}
\wp(z, \tau) = \pi^2 \theta_2^2(0,\tau) \theta_3^2(0,\tau)
    \frac{\theta_4^2(z,\tau)}{\theta_1^2(z,\tau)} -
    \frac{\pi^2}{3} \left[ \theta_2^4(0,\tau) + \theta_3^4(0,\tau)\right].
\label{eq:weierp}
\end{equation}
For all these functions, Arb also allows computing an arbitrary number of $z$-derivatives.
Derivatives are handled by working with $\theta_j(z+x,\tau)$ and $\wp(z+x,\tau)$
as elements of $\mathbb{C}[[x]]$ (truncated to some length $O(x^D)$),
using power series arithmetic.

Arb also implements the quasielliptic Weierstrass zeta and sigma functions $\zeta(z,\tau)$ and $\sigma(z,\tau)$
as well as the lattice invariants $g_2, g_3$
(which are essentially Eisenstein series)
and lattice roots $4z^3 - g_2 z - g_3 = 4(z-e_1)(z-e_2)(z-e_3)$
arising in the differential equation $[\wp'(z, \tau)]^2 = 4 [\wp(z,\tau)]^3 - g_2 \wp(z,\tau) - g_3$.
The inverse Weierstrass elliptic function is also available; see section~\ref{sect:incellip}.

\index{Eisenstein series}
\index{lattice invariant}
\index{lattice root}
\index{Jacobi elliptic function}
\index{differential equation}

The Jacobi elliptic functions sn, cn, $\ldots$
are not currently part of the library,
but users can compute them via theta functions using formulas similar to \eqref{eq:weierp}.

\subsection{Argument reduction}

\label{sect:thetatransform}

As the first step when computing theta functions or elliptic functions,
we reduce~$\tau$ to the fundamental domain using modular transformations.
This gives us a new lattice parameter $\tau'$
and a new argument $z'$.
As a second step, we reduce $z'$ modulo $\tau'$,
giving an argument $z''$ with smaller imaginary part (it is not necessary to reduce $z'$ modulo 1
since this is captured by the oscillatory part of exponentials).
We can then compute $\theta_j(z'', \tau')$ using the theta series \eqref{eq:genthetaseries}.

\index{argument reduction}
\index{lattice parameter}
\index{elliptic function}
\index{Jacobi theta function}
\index{lattice}
\index{modular transformation}

These steps together ensure that both $|q|$ and $\max(|w|,|w|^{-1})$ will be small.
It is important to perform both transformations.
Consider $\tau = 0.07 + 0.003i$ and $z = 3.14+2.78i$:
without the modular transformation, the direct series evaluation
would use 3710 terms for machine precision.\footnote{The number is somewhat smaller if the series is truncated optimally using a relative rather than an absolute tolerance.}
With the modular transformation alone, it would use 249 terms.
With both reductions, only 6 terms are used!
Depending on the arguments, the numerical stability may also be
improved substantially.

Modular transformations have the effect of permuting
the theta functions and introducing certain exponential prefactors.
It is easy to write down the transformations for
the generators $\tau+1$, $-1/\tau$,
but the action of a composite transformation involves a certain
amount of bookkeeping.
The steps have been worked out by Rademacher~\cite[chapter 10]{Rademacher1973}.
We reproduce the formulas below.%
\footnote{We give the inverse form of the transformation.}

We wish to write a theta function with lattice parameter $\tau$ in terms
of a theta function with lattice parameter $\tau' = g(\tau)$, given
some $g = \smash{\begin{psmallmatrix} a & b \\ c & d \end{psmallmatrix}} \in \text{PSL}(2, \mathbb{Z})$.
For $j = 0, 1, 2, 3$, there are $R_j, S_j \in \mathbb{Z}$ depending on $g$ such that
\begin{equation}
\label{eq:thetatransform}
\theta_{1+j}(z,\tau) = \exp(\pi i R_j / 4) \cdot A \cdot B \cdot \theta_{1+S_j}(z',\tau')
\end{equation}
where if $c = 0$,
\begin{equation}
z' = z, \quad A = 1, \quad B = 1, \\
\end{equation}
and otherwise (if $c > 0$),
\begin{equation}
z' = \frac{-z}{c \tau + d}, \quad
A = \sqrt{\frac{i}{c \tau + d}}, \quad
B = \exp\left(-\pi i c \frac{z^2}{c \tau + d}\right).
\end{equation}
We always have $B = 1$ when
computing theta constants which have $z = 0$.

The parameters $R_j, S_j$ are computed from $g$ as follows.
If $c = 0$, we have
$\theta_j(z, \tau) = \exp(-\pi i b / 4) \theta_j(z, \tau+b)$
for $j = 1, 2$, whereas
$\theta_3$ and $\theta_4$ remain unchanged when $b$ is even
and swap places with each other when $b$ is odd.
For the $c > 0$ case, it is helpful to define the function $\theta_{m,n}(z,\tau)$ for $m, n \in \mathbb{Z}$ by
\begin{align}
\theta_{0,0}(z,\tau) &= \theta_3(z,\tau), & \theta_{0,1}(z,\tau) &= \theta_4(z,\tau), \\ \nonumber
\theta_{1,0}(z,\tau) &= \theta_2(z,\tau), & \theta_{1,1}(z,\tau) &= i \theta_1(z,\tau), \\
\theta_{m+2,n}(z,\tau) &= (-1)^n \theta_{m,n}(z,\tau) & \theta_{m,n+2}(z,\tau) &= \theta_{m,n}(z,\tau). \nonumber
\end{align}
With this notation, we have
\begin{align}
 \theta_1(z,\tau) &= \varepsilon_1 A B \theta_1(z', \tau'), &
 \theta_2(z,\tau) &= \varepsilon_2 A B \theta_{1-c,1+a}(z', \tau'), \\ \nonumber
 \theta_3(z,\tau) &= \varepsilon_3 A B \theta_{1+d-c,1-b+a}(z', \tau'), &
 \theta_4(z,\tau) &= \varepsilon_4 A B \theta_{1+d,1-b}(z', \tau')
\end{align}
where $\varepsilon_k$ is an 8th root of unity.
If we denote by $\varepsilon(a,b,c,d) = \exp(\pi i R(a,b,c,d) / 12)$
the 24th root of unity
in the transformation \eqref{eq:etatransform} of the Dedekind eta function, then
\begin{align} \nonumber
  \varepsilon_1(a,b,c,d) &= \exp(\pi i [R(-d,b,c,-a) + 1] / 4), \\ \nonumber
  \varepsilon_2(a,b,c,d) &= \exp(\pi i [-R(a,b,c,d) + (5+(2-c)a)] / 4), \\
  \varepsilon_3(a,b,c,d) &= \exp(\pi i [-R(a,b,c,d) + (4+(c-d-2)(b-a))] / 4), \\ \nonumber
  \varepsilon_4(a,b,c,d) &= \exp(\pi i [-R(a,b,c,d) + (3-(2+d)b)] / 4).  \\ \nonumber
\end{align}
Finally, to reduce $z'$, we compute
$n = \lfloor \operatorname{Im}(z') / \operatorname{Im}(\tau') + 1/2 \rfloor$
and set $z'' = z' - n \tau'$.
In this step, all theta functions pick up a prefactor
$\exp(\pi i [-\tau n^2 - 2nz])$ (this data may be combined with $B$)
while $\theta_1$ and $\theta_2$ pick up the additional prefactors $(-1)^n$
(this data may be combined with $R_j$).

\index{Dedekind eta function}
\index{power series arithmetic}
\index{modular transformation}
\index{ball arithmetic}
\index{Arb library}

When computing $z$-derivatives of theta functions, the same formulas are applied
in power series arithmetic.
That is, if the initial argument consists of the formal power series
$z + x$, then the scaling factor $-1/(c \tau + d)$ is applied
coefficient by coefficient, while
$B = B_0 + B_1 x + \ldots$ is obtained by squaring the power series $z + x$, scaling,
and then evaluating a power series exponential.

As with modular forms, the transformations should be applied at the highest possible
level. For example, when computing a quotient of two theta functions
of the same $z, \tau$, the prefactors $A$ and $B$ in \eqref{eq:thetatransform}
cancel out (and the leading roots of unity possibly also simplify).
We should then simplify the expression symbolically and
avoid computing $A$ and $B$ altogether, since this both saves time
and improves numerical stability in ball arithmetic (in particular,
$e^{f(z)} / e^{f(z)}$ evaluated in ball arithmetic will not give~1 but rather
a ball which can be extremely wide).

Since the description of the algorithm given above is quite terse,
the reader may find it helpful to look at the code in Arb to see the concrete
steps.

\subsection{Theta function series evaluation}

\index{Jacobi theta function}

Algorithm~\ref{alg:theta} implements the expansions~\eqref{eq:genthetaseries},
with the optimization that we combine operations to save work
when computing all four functions and their derivatives simultaneously
(a single theta function could be computed slightly faster,
but computing all four functions is barely more work
than it would be to compute a pair containing either $\theta_1$ or $\theta_2$ and either $\theta_3$ or $\theta_4$).
This is essentially the algorithm used in Arb for $z \ne 0$,
while more optimized code is used for theta constants.

The main index $k$ runs over the terms in the following order:
$$
\arraycolsep=10pt
\begin{array}{llll}
    & \theta_1, \theta_2 & q^0 & (w^1 \pm w^{-1}) \\
k = 0  & \theta_3, \theta_4 & q^1 & (w^2 \pm w^{-2}) \\
k = 1  & \theta_1, \theta_2 & q^2 & (w^3 \pm w^{-3}) \\
k = 2  & \theta_3, \theta_4 & q^4 & (w^4 \pm w^{-4}) \\
k = 3  & \theta_1, \theta_2 & q^6 & (w^5 \pm w^{-5}) \\
k = 4  & \theta_3, \theta_4 & q^9 & (w^6 \pm w^{-6}) \\
k = 5  & \theta_1, \theta_2 & q^{12} & (w^7 \pm w^{-7}) \\
\end{array}
$$

The algorithm outputs the range of scaled derivatives $\theta_j^{(r)}(z,\tau) / r!$ for $0 \le r < D$.
The term of index $k$ in the main summation picks up a factor
$\pm (k+2)^r$ from $r$-fold differentiation of $w^{k+2}$.
Another factor $(\pi i)^r /r!$ is needed to convert to
a $z$-derivative and a power series coefficient,
but we postpone this to a single rescaling pass at the end of the computation.
In the main summation, we write the even
cosine terms as $w^{2n} + w^{-2n}$, the odd cosine terms as
$w (w^{2n} + w^{-2n-2})$, and the sine terms as $w (w^{2n} - w^{-2n-2})$,
postponing a multiplication by $w$ for $\theta_1$ and $\theta_2$ until the end,
so that only even powers of $w$ and $w^{-1}$ are needed.

For some integer $N \ge 1$, the summation is stopped just before term
$k = N$. Let $Q = |q|$, $W = \max(|w|,|w^{-1}|)$,
$E = \lfloor (N+2)^2 / 4 \rfloor$ and 
$F = \lfloor (N+1)/2 \rfloor + 1$. The error of the
zeroth derivative can be bounded as
\begin{equation}
        2 Q^E W^{N+2} \left[ 1 + Q^F W + Q^{2F} W^2 + \ldots \right]
        = \frac{2 Q^E W^{N+2}}{1 - Q^F W}
\end{equation}
    provided that the denominator $1 - Q^F W$ is positive.
For the $r$-th derivative, including the factor $(k+2)^r$
gives the error bound
\begin{equation}
        2 Q^E W^{N+2} (N+2)^r \left[ 1 + Q^F W \frac{(N+3)^r}{(N+2)^r} + Q^{2F} W^2 \frac{(N+4)^r}{(N+2)^r} + \ldots \right]
\end{equation}
    which by the inequality $(1 + m/(N+2))^r \le \exp(mr/(N+2))$
    can be bounded as
\begin{equation}
\frac{2 Q^E W^{N+2} (N+2)^r}{1 - Q^F W \exp(r/(N+2))},
\end{equation}
again valid when the denominator is positive.

\pagebreak

\begin{algorithm}[h!]
  \caption{Computation of Jacobi theta functions (using series evaluation)}
  \small
  \label{alg:theta}
  \begin{algorithmic}[1]
  \Require $z, \tau \in \mathbb{C}$ with $\operatorname{Im}(\tau) > 0$ (can be arbitrary, but should be reduced for best performance), integer $D \ge 1$ to output the $D$ first terms in the Taylor expansions with respect to $z$, precision $p$
  \Ensure $\theta_j = [\alpha_0, \ldots, \alpha_{D-1}]$ represents $\theta_j(z+x,\tau) = \alpha_0 + \alpha_1 x + \ldots + \alpha_{D-1} x^{D-1}$, for $1 \le j \le 4$
  \State $q_4 \gets e^{\pi i \tau/4}; \; q \gets q_4^4; \; w \gets e^{\pi i z}; \; v \gets w^{-1}; \; Q \gets |q|; \; W \gets \max(|w|,|v|)$
  \State Choose $N$ with $E = \lfloor (N+2)^2 / 4 \rfloor$ and $F = \lfloor (N+1)/2 \rfloor + 1$ such that $Q^E W^{N+2} < 2^{-p}$ and $\alpha = Q^F W \exp(r/(N+2)) < 1$
  \State $\textbf{ for } 0 \le r < D$ $\textbf{ do }$ $\varepsilon[r] \gets 2 Q^E W^{N+2} (N+2)^r / (1 - \alpha)$ $\textbf{ end for }$ \Comment{Error bounds}
  \State $\textbf{w} \gets [1, w^2, w^4, \ldots, w^{2K-2}]; \; \textbf{v} \gets [1, v^2, v^4, \ldots, v^{2K}]$ for $K = \lfloor (N + 3) / 2 \rfloor$ \Comment{Precompute powers}
  \State $\theta_1 \gets [0,\ldots,0]; \; \theta_2 \gets [0,\ldots,0]; \; \theta_3 \gets [0,\ldots,0]; \; \theta_4 \gets [0,\ldots,0]; \;$ \Comment{Arrays of length $D$}
  \For {$0 \le k < N$}
    \State $m \gets \lfloor (k+2)^2 / 4 \rfloor; \; n \gets \lfloor k / 2 \rfloor + 1$
    \State Compute $q^m$ \Comment{Use addition sequence~\cite[Alg. 2]{enge2018short} to build $q^m$ from previous powers.}
    \State $t \gets (\textbf{w}[n] + \textbf{v}[n + (k \bmod 2)]) q^m$
    \State $u \gets (\textbf{w}[n] - \textbf{v}[n + (k \bmod 2)]) q^m$ \Comment{Skip when $k \bmod 2 = 0$ if $D = 1$.}
    \If {$k \bmod 2 = 0$}
        \For {$0 \le r < D$}
          \If {$r \bmod 2 = 0$}
            \State $\textbf{ if } r \ne 0 \textbf{ then } t \gets 4 n^2 t \textbf{ end if}$
            \State $\theta_3[r] \gets \theta_3[r] + t; \; \theta_4[r] \gets \theta_4[r] + (-1)^{\lfloor (k+2)/2 \rfloor } t$
          \Else
            \State $\textbf{ if } r = 1 \textbf{ then } u \gets 2 n u \textbf{ else } u \gets 4 n^2 u \textbf{ end if}$
            \State $\theta_3[r] \gets \theta_3[r] + u; \; \theta_4[r] \gets \theta_4[r] + (-1)^{\lfloor (k+2)/2 \rfloor } u$
          \EndIf
        \EndFor
    \Else
        \For {$0 \le r < D$}
          \If {$r \bmod 2 = 0$}
            \State $\theta_1[r] \gets \theta_1[r] + (-1)^{\lfloor (k+1)/2 \rfloor } u; \; \theta_2[r] \gets \theta_2[r] + t$
            
          \Else
            \State $\theta_1[r] \gets \theta_1[r] + (-1)^{\lfloor (k+1)/2 \rfloor } t; \; \theta_2[r] \gets \theta_2[r] + u$
          \EndIf
          \State $t \gets (2n+1) t; \; u \gets (2n+1) u$
        \EndFor
    \EndIf
  \EndFor
  \For {$0 \le r < D$}
    \State $\theta_1[r] \gets \theta_1[r] w + (w - (-1)^r v)$ \Comment{Adjust power of $w$ and add leading terms}
    \State $\theta_2[r] \gets \theta_2[r] w + (w + (-1)^r v)$
    \State $\textbf{ for } 1 \le j \le 4$ $\textbf{ do }$ $\theta_j[r] \gets \theta_j[r] + [\pm \varepsilon[r]] + [\pm \varepsilon[r]] i$ $\textbf{ end for }$ \Comment{Add error bounds}
    \State $C \gets (\pi i)^r / r!$ \Comment{Final scaling factors}
    \State $\theta_1[r] \gets -i q_4 C \theta_1[r]; \; \theta_2[r] \gets q_4 C \theta_2[r]; \; \theta_3[r] \gets C \theta_3[r]; \; \theta_4[r] \gets C \theta_4[r]$
  \EndFor
  \State $\theta_3[0] \gets \theta_3[0] + 1; \; \theta_4[0] \gets \theta_4[0] + 1$ \Comment{Add leading terms}
  \end{algorithmic}
\end{algorithm}

The time complexity of the algorithm is $\OOtilde(p^{1.5})$ (with all inputs besides $p$ fixed).
By employing fast Fourier transforms cleverly, \index{Fast Fourier transform}
the complexity of evaluating theta functions from their series expansions
can be reduced to $\OOtilde(p^{1.25})$, but
that method is only faster in practice for $p$ exceeding
200\,000 bits~\cite{nogneng2018evaluation}. See also section~\ref{sect:quadratic}
below concerning methods that are even faster asymptotically.

\section{Complete elliptic integrals and the AGM}
\label{sect:complagm}
\vspace{1mm}\noindent
Complete elliptic integrals arise in period relations
for elliptic functions.
The complete elliptic integral of the first kind
is
$K(m) = \tfrac{1}{2} \pi {}_2F_1(\tfrac{1}{2}, \tfrac{1}{2}, 1, m)$ and
the complete elliptic integral of the second kind is
$E(m) = \tfrac{1}{2} \pi {}_2F_1(-\tfrac{1}{2}, \tfrac{1}{2}, 1, m)$
where ${}_2F_1$ denotes the Gauss hypergeometric function,
defined for $|z| < 1$ by
\begin{equation}
{}_2F_1(a,b,c,z) = \sum_{k=0}^{\infty} \frac{(a)_k (b)_k}{(c)_k} \frac{z^k}{k!}, \quad (x)_k = x (x+1) \cdots (x+k-1)
\end{equation}
and elsewhere by analytic continuation with the standard branch cut on $[1,\infty)$.

\index{elliptic integral}
\index{complete elliptic integral}
\index{analytic continuation}
\index{hypergeometric series}
\index{Gauss hypergeometric function}

The ${}_2F_1$ function
can be computed efficiently for any $z \in \mathbb{C}$ using a combination
of argument transformations, analytic continuation techniques,
and series expansions (where the rectangular splitting trick~\eqref{eq:rectsplit}
and other accelerations methods are applicable).
A general implementation of ${}_2F_1$
exists in Arb~\cite{Johansson2016hypergeometric}.
However, it is more efficient to compute
the complete elliptic
integrals by exploiting their connection with the
arithmetic-geometric mean (AGM) described below.

\index{arithmetic-geometric mean}

A third complete elliptic integral $\Pi(n,m)$ is also
encountered, but this is a more
complicated function that is not a special case of ${}_2F_1$,
and we handle it later in terms of an incomplete integral
without using a dedicated algorithm for the complete case.

The arithmetic-geometric mean $M(x,y)$ of two nonnegative real numbers $x, y$
is defined as the common limit of the sequences
\begin{equation}
a_{n+1} = \frac{a_n+b_n}{2}, \quad b_{n+1} = \sqrt{a_n b_n}
\end{equation}
with initial values $a_0 = x, b_0 = y$. In different
words, the AGM can be
computed by repeatedly applying
the functional equation
$M(x,y) = M\left((x+y)/2, \sqrt{x y}\right)$.
It is a well known fact that each step of the AGM iteration
roughly doubles the number of accurate digits in the approximation $a_n \approx b_n \approx M(x,y)$,
so it only costs $O(\log(p))$ arithmetic operations
to compute the AGM with an accuracy of $p$ bits,
resulting in a bit complexity of $\OOtilde(p)$.

For complex $x$, $y$, defining the AGM becomes
more difficult since there are two possible choices for the square root
in each step of the iteration, and these choices lead to different limits.
However, it turns out that there is an ``optimal'' choice which leads
to a well-defined and useful extension of the AGM to complex variables.
We rely on several properties of this function
proved in earlier work~\cite{cox2000arithmetic,dupont2006moyenne,cremona2013complex}.

With complex variables, it is convenient to work with the univariate
function $M(z) = M(1,z)$, with a
branch cut on $(-\infty,0]$. The general case can be recovered as
$M(x,y) = x M(1,y/x)$.
The complete elliptic integrals (with the conventional
branch cuts on $[1,\infty)$)
are now given by
\begin{equation}
K(m) = \frac{\pi}{2 M(\sqrt{1-m})}, \quad E(m) = (1-m) (2mK'(m) + K(m)).
\end{equation}

For implementing the function $M(z)$, we can further assume that $\operatorname{Re}(z) \ge 0$ holds.
If this is not the case, we first apply
the functional equation $M(z) = (z+1) M(u) / 2$
where $u = \sqrt{z} / (z+1)$.
The correct square root in the AGM iteration
is now always equal to $\sqrt{a_n} \sqrt{b_n}$,
written in terms of the usual principal square root function.
This can be computed as $\sqrt{a_n b_n}$, $i \sqrt{-a_n b_n}$, $-i \sqrt{-a_n b_n}$, $\sqrt{a_n} \sqrt{b_n}$
respectively if both $a_n$ and $b_n$ have positive real part, nonnegative
imaginary part, nonpositive imaginary part, or otherwise.
When the iteration is executed in ball arithmetic,
the computed balls may end up containing
points with negative real part,
but this just inflates the final result and does not affect correctness.

\index{argument reduction}

The iteration should be terminated when $a_n$ and $b_n$ are close enough.
For positive real variables, we can simply take lower and upper bounds
to get a correct enclosure. For complex variables, it can be shown
\cite[p.~87]{dupont2006moyenne} that
$|M(z) - a_n| \le |a_n - b_n|$ if $\operatorname{Re}(z) \ge 0$, giving a convenient error bound.
However, instead of running the AGM iteration until $a_n$ and $b_n$ agree to
$p$ bits, it is slightly better to stop when
they agree to about $p/10$ bits and end with a Taylor series.
With $t = (a-b)/(a+b)$, we have
\begin{equation}
M(a,b) = \frac{(a+b) \pi}{4 K(t^2)}, \quad \frac{\pi}{4 K(t^2)} = \tfrac{1}{2} - \tfrac{1}{8} t^2 - \tfrac{5}{128} t^4 - \tfrac{11}{512} t^6 - \tfrac{469}{32768} t^8 + \ldots
\end{equation}
which is valid at least when $|t| < 1$ and $a, b$ have nonnegative real part,
and where the tail ($\ldots$) is bounded by $\sum_{k=10}^{\infty} |t|^k/64$.

\index{quadratic convergence}

This algorithm follows the
pattern of argument reduction and series evaluation.
However, unlike the elementary functions and the incomplete
elliptic integrals described below,
there is no asymptotic benefit to using more terms of the series.
The quadratic convergence of the AGM iteration is so rapid
that we only get a speedup from trading $O(1)$ of the $O(\log p)$ square roots
for lower-overhead multiplications.
Although there is no asymptotic improvement,
the order-10 series expansion nevertheless
gives a significant speedup up to a few thousand bits.

For computing the second elliptic integral
$E(m)$ or the first derivative $M'(z)$ of the AGM,
a simple method is to use a central finite difference to compute
$(M(z), M'(z)) \approx (M(z+h)+M(z-h))/2, (M(z+h)-M(z-h))/(2h)$.
This requires two evaluations at 1.5 times increased precision,
which is about three times as expensive as evaluating $M$ once.
Error bounds can be obtained using the Cauchy integral formula
and the inequality $|M(z)| \le \max(1, |z|)$ which is
an immediate consequence of the AGM iteration.
This method has been implemented in Arb.
A more efficient method is to compute $E(m)$ using an auxiliary
sequence related to the AGM iteration,
which also generalizes to computing $\Pi(n,m)$~\cite[19.8.6~and~19.8.7]{Olver2010}.
This method has not yet been implemented in Arb
since it requires some additional error analysis
and study for complex variables.

Higher derivatives of the arithmetic-geometric mean
or the complete elliptic integrals can be computed using recurrence relations.
Writing $W(z) = 1 / M(z)$
and $W(z+x) = \sum_{k=0}^{\infty} c_k x^k$,
we have
$-2 z (z^2-1) c_2 = (3z^2-1) c_1 + z c_0$,
$-(k+2)(k+3) z (z^2-1) c_{k+3} = (k+2)^2 (3z^2-1) c_{k+2} + (3k(k+3)+7)z c_{k+1} + (k+1)^2 c_{k}$
when $z \ne 1$ and
$-(k+2)^2 c_{k+2} = (3k(k+3)+7) c_{k+1} + (k+1)^2 c_{k}$
when $z = 1$.

\section{Incomplete elliptic integrals}
\label{sect:incellip}

\vspace{1mm}\noindent
A general elliptic integral is an integral of the
form $\int_a^b R(t,\sqrt{P(t)}) dt$ where $R$ is a bivariate
rational function and $P$ is a cubic or quartic polynomial
without repeated roots. It is well known that
any elliptic integral can be expressed
in terms of integrals of rational functions and a finite set of
standard elliptic integrals.

\index{elliptic integral}
\index{incomplete elliptic integral}
\index{Legendre form}

Such a set of standard integrals is given by the Legendre incomplete elliptic integrals
of the first, second and third kind
\begin{equation}
F(\phi,m) = \int_0^{\phi} \frac{dt}{\sqrt{1-m \sin^2 t}}, \quad E(\phi,m) = \int_0^{\phi} \sqrt{1-m \sin^2 t} \, dt,
\label{eq:legendreinc1}
\end{equation}
\begin{equation}
\Pi(n, \phi, m) = \int_0^{\phi}
    \frac{dt}{(1-n \sin^2 t) \sqrt{1-m \sin^2 t}}.
\label{eq:legendreinc2}
\end{equation}

The complete elliptic integrals are the special cases $E(m) = E(\pi/2,m)$,
$K(m) = F(\pi/2,m)$, and $\Pi(n,m) = \Pi(n,\pi/2,m)$.

The definitions for complex variables do not appear to be standardized
in the literature,
but following the conventions used in Mathematica~\cite{wolfellip}, we may
fix an interpretation of \eqref{eq:legendreinc1}--\eqref{eq:legendreinc2}
on $-\pi/2 \le \operatorname{Re}(\phi) \le \pi/2$ and use
the quasiperiodic extensions
$F(\phi + k \pi, m) = 2 k K(m) + F(\phi,m)$,
$E(\phi + k \pi, m) = 2 k E(m) + E(\phi,m)$,
$\Pi(n, \phi + k \pi, m) = 2 k \Pi(n,m) + \Pi(n,\phi,m)$
for $k \in \mathbb{Z}$.\footnote{For $\Pi$, Mathematica
restricts this quasiperiodicity relation to hold
only for $-1 \le n \le 1$.}

The Legendre forms of incomplete elliptic integrals are widely
used by tradition, but they have some practical drawbacks.
Since they have a complicated (and not standardized) complex branch
structure, transforming their arguments using functional equations
or using them to represent other functions often requires making
complicated case distinctions.
As a result, it is cumbersome both to compute the functions themselves
and to apply them, outside of a restricted parameter range.

\index{analytic continuation}
\index{Appell hypergeometric function}
\index{Lauricella hypergeometric function}
\index{hypergeometric series}

We remark that $F$ and $E$ can be expressed in terms of the
Appell $F_1$ hypergeometric function of two variables, while
$\Pi$ can be expressed in terms of the three-variable Lauricella
hypergeometric function \smash{$F_D^{(3)}$}, generalizing the ${}_2F_1$
representations for the complete integrals. Such formulas are by themselves mainly useful
when the hypergeometric series converge,
and provide no insight into the analytic continuations.

\index{Carlson symmetric form}

In the 1960s, Carlson introduced an alternative set of
standard elliptic integrals in which all or some of the variables are symmetric~\cite{carlson1995numerical}.
The Carlson incomplete elliptic integrals are
\begin{equation}
R_F(x,y,z) = \frac{1}{2} \int_0^{\infty} \frac{dt}{\sqrt{(t+x)(t+y)(t+z)}}
\label{eq:carlson1}
\end{equation}
and
\begin{equation}
R_J(x,y,z,p) = \frac{3}{2} \int_0^{\infty} \frac{dt}{(t+p)\sqrt{(t+x)(t+y)(t+z)}}
\label{eq:carlson2}
\end{equation}
together with three special cases $R_D(x, y, z) = R_J(x, y, z, z)$,
$R_C(x, y) = R_F(x, y, y)$, and
\begin{equation}
R_G(x,y,z) = z R_F(x,y,z) - \frac{1}{3} (x-z)(y-z) R_D(x,y,z) + \frac{\sqrt{x} \sqrt{y}}{\sqrt{z}}.
\end{equation}

The Carlson forms have several advantages over the Legendre forms.
Symmetry unifies and simplifies the argument transformation formulas,
and the Carlson forms also have a simpler complex branch structure,
induced by choosing the branch of the square root
in \eqref{eq:carlson1} and \eqref{eq:carlson2}
to extend continuously from $+\infty$.
We can define and compute the Legendre forms from the Carlson forms using
\begin{align} \nonumber
F(\phi,m) &= s R_F(x, y, 1), \\
E(\phi,m) &= s R_F(x,y,1) - \tfrac{1}{3} m s^3 R_D(x,y,1), \\
\Pi(n,\phi,m) &= s R_F(x,y,1) + \tfrac{1}{3} n s^3 R_J(x,y,1,p) \nonumber
\end{align}
on $-\pi/2 \le \operatorname{Re}(\phi) \le \pi/2$ (with the quasiperiodic
extensions elsewhere)
where $x = c^2$, $y=1-m s^2$, $p = 1-n s^2$ and $s = \sin(\phi)$, $c = \cos(\phi)$.
This is the approach used to implement the Legendre forms in Arb.
The Carlson forms themselves are also exposed to users.
Formulas for other elliptic integrals can be found in~\cite{carlson1995numerical}.

\index{Arb library}
\index{incomplete elliptic integral}
\index{elliptic integral}
\index{Carlson symmetric form}
\index{Weierstrass elliptic function}

Elliptic integrals can also be characterized as the inverse
functions of elliptic functions.
For example, the inverse of the Weierstrass elliptic function,
which by definition satisfies $\wp(\wp^{-1}(z,\tau),\tau) = z$,
is given by the elliptic integral
\begin{equation}
\wp^{-1}(z, \tau) = \frac{1}{2} \int_z^{\infty} \frac{dt}{\sqrt{(t-e_1)(t-e_2)(t-e_3)}}
    = R_F(z-e_1,z-e_2,z-e_3).
\end{equation}
The implementation in Arb simply computes
the lattice roots $e_1,e_2,e_3$ using theta constants and then calls $R_F$.
The inverses of Jacobi's elliptic functions can be computed
similarly, but at this time they are not implemented in Arb.

\index{Weierstrass elliptic function}
\index{Jacobi elliptic function}
\index{lattice root}

Carlson gives algorithms for computing $R_F$ and $R_J$
using argument reduction and series evaluation~\cite{carlson1995numerical}.
The algorithm for $R_F$ is correct for all complex $x,y,z$
(Carlson restricts to the cut plane with $(-\infty,0)$ removed, but
it is clear that the algorithm also works on the branch cut by continuity).
The algorithm for $R_J$ is not correct for all values of the variables,
but it is always correct when computing $R_D$
(otherwise, a sufficient condition is that $x,y,z$ have nonnegative real part
while $p$ has positive real part).
Carlson also provides modifications
of the algorithms for computing the Cauchy principal values of the integrals.

\index{ball arithmetic}
\index{Arb library}

We will now describe Carlson's algorithm for $R_F$ and adapt
it to the setting of arbitrary precision and ball arithmetic.
The algorithm given in~\cite{carlson1995numerical} for $R_J$ and $R_D$
works analogously, but we do not reproduce all the steps here
since the formulas would be too lengthy
(the code in Arb can be consulted for
concrete details).

\subsection{Argument reduction}

Argument reduction for $R_F$ uses the symmetric ``duplication formula''
\begin{equation}
R_F(x,y,z) = R_F\left(\frac{x+\lambda}{4}, \frac{y+\lambda}{4}, \frac{z+\lambda}{4}\right)
\label{eq:rfdup}
\end{equation}
where $\lambda = \sqrt{x} \sqrt{y} + \sqrt{y} \sqrt{z} + \sqrt{z} \sqrt{x}$.
Each application of \eqref{eq:rfdup} reduces the distance between the
arguments by roughly a factor 4.
The analogous formula for $R_J$ reads
\begin{equation}
R_J(x,y,z,p) = \frac{1}{4} R_J\left(\frac{x+\lambda}{4}, \frac{y+\lambda}{4}, \frac{z+\lambda}{4}, \frac{p+\lambda}{4}\right)
+ \frac{6}{d} \, R_C(1,1+e)
\label{eq:rjdup}
\end{equation}
where $\lambda$ is defined as above
and $d, e$ are certain auxiliary terms (see~\cite[(24)--(28)]{carlson1995numerical}).
The formulas \eqref{eq:rfdup} and \eqref{eq:rjdup}
are iterated until all the parameters are close
so that a series expansion can be used, as detailed in the next subsection.
It is interesting to note the similarity between \eqref{eq:rfdup}
and the AGM iteration, although the convergence rate of \eqref{eq:rfdup}
only is linear.

\index{argument reduction}
\index{Carlson symmetric form}
\index{arithmetic-geometric mean}

When computing $R_C$ or $R_D$, some redundant operations
in the reductions for $R_F$ and $R_J$ can be avoided.
$R_C(x,y)$ can also be expressed piecewise using inverse
trigonometric and hyperbolic functions.
The special case
$R_C(1, 1+t) = \operatorname{atan}(\sqrt{t})/\sqrt{t} = {}_2F_1(1,\tfrac{1}{2},\tfrac{3}{2},-t)$
is particularly important, as it
is needed in the evaluation of $R_J$.
This function is better computed via the inverse tangent function
(or a direct Taylor series for small $|t|$)
than by invoking Carlson's general method for $R_F$.

\subsection{Series expansions}

\index{hypergeometric series}

Carlson's incomplete elliptic integrals are special cases of
a multivariate hypergeometric function that may be written as
\begin{equation}
R_{-a}(z_1,\ldots,z_n) = A^{-a} \sum_{N=0}^{\infty} \frac{(a)_n}{(\tfrac{1}{2}n)_N} T_N(Z_1,\ldots,Z_n)
\label{eq:raseries}
\end{equation}
where $A = \tfrac{1}{n} \sum_{j=1}^n z_j$, $Z_j = 1 - z_j / A$, and
\begin{equation}
T_N(Z_1,\ldots,Z_n) = \sum_{\substack{m_1+\ldots+m_n=N \\ m_1,\ldots,m_n \ge 0}} \; \prod_{j=1}^n \frac{(\tfrac{1}{2})_{m_j}}{m_j!} Z_j^{m_j}.
\label{eq:raterm}
\end{equation}

We have $R_F(x,y,z) = R_{-1/2}(x,y,z)$ and $R_J(x,y,z,p) = R_{-3/2}(x,y,z,p,p)$.
The crucial property of this hypergeometric representation is that the
expansion point is the arithmetic mean of the arguments.
After sufficiently many argument reduction steps have been performed,
we will have $z_1 \approx z_2 \approx \ldots \approx z_n$
which ensures $|Z_1|,\ldots,|Z_n| \ll 1$ and rapid convergence of the series.
A trivial bound for the terms is
$|T_N(Z_1,\ldots,Z_n)| \le p(N) \max(|Z_1|,\ldots,|Z_n|)^N$,
where $p(N)$ denotes the number of partitions of $N$ which
is bounded by $O(c^N)$ for any $c > 1$.
An explicit calculation shows, for example, that the error
when computing either $R_F$ or $R_J$ is bounded by
\begin{equation}
2 A^{-a} \sum_{N=B}^{\infty} \left(\tfrac{9}{8} \max(|Z_1|,\ldots,|Z_n|)\right)^N
\end{equation}
if the summation in \eqref{eq:raseries} includes the terms of order $N < B$.

\index{elementary symmetric polynomial}
\index{hypergeometric series}

For the evaluation of \eqref{eq:raterm}, Carlson noted that it is
more efficient to work with elementary symmetric polynomials $E_j = E_j(Z_1,\ldots,Z_n)$
instead of the direct variables $Z_1, \ldots, Z_n$, giving
\begin{equation}
T_N(Z_1,\ldots,Z_n) = \sum_{\substack{m_1+2m_2+\ldots+n m_n=N \\ m_1,\ldots,m_n \ge 0}} \; (-1)^{M+N} (\tfrac{1}{2})_M \prod_{j=1}^n \frac{E_j^{m_j}}{m_j!}.
\label{eq:rasymterm}
\end{equation}
The key observation is that the symmetric choice of expansion variables $Z_j$
with respect to $z_1,\ldots,z_n$ implies that $E_1 = 0$, which eliminates
most of the terms in \eqref{eq:rasymterm}.\footnote{This is an algebraic simplification, so we can take $E_1 = 0$ even if the
input argument are represented by inexact balls.}
This dramatically reduces the amount of work to compute $T_N$
compared to \eqref{eq:raterm}.
For the $R_{-1/2}$ series, there are
$(N+1) (N+2) / 2$ terms in \eqref{eq:raterm} and roughly $N / 6$ terms in \eqref{eq:rasymterm};
for example, if $N = 8$, there are 45 terms in the former
and only two nonzero terms (with monomials $E_1 E_2^2$ and $E_1^4$) in the latter.

\subsection{Series evaluation and balanced argument reduction}
\label{sect:rfalg}

The argument reduction effectively adds $2B$ bits per step
if a series expansion of order $B$ is used.
In other words, roughly $p / (2B)$ argument reduction steps are needed
for $p$-bit precision.
Carlson suggests using a precomputed truncated series
of order $B = 6$ or $B = 8$,
which is a good default at machine precision and up to a few hundred bits.
At higher precision, we make the observation that
it pays off to vary $B$ dynamically as a function of $p$
and evaluate the series with an algorithm.

\begin{algorithm}[h!]
  \caption{Computation of $R_F(x,y,z)$}
  \small
  \label{alg:hyprs}
  \begin{algorithmic}[1]
    \State Choose series truncation order $B$ optimally depending on the precision $p$
    \State Apply argument reduction \eqref{eq:rfdup} until $x,y,z$ are close (until $\varepsilon \approx 2^{-p}$ below)
    \State $A \gets (x+y+z)/3; \;\; (X,Y,Z) \gets (1-x/A, 1-y/A, 1-z/A); \;\; (E_2, E_3) \gets (XY-Z^2, XYZ)$
    \State $\varepsilon \gets 2 \sum_{k=B}^{\infty} \left(\tfrac{9}{8} \max(|X|,|Y|,|Z|)\right)^k$ \Comment{Series error bound}
    \Procedure{RSum}{$E_2,E_3,B$} \Comment{Compute $R = \sum_{N=0}^{B-1} \left[ (\tfrac{1}{2})_N / (\tfrac{3}{2})_N \right] T_N$}
    \State Precompute $E_2^k$ for $2 \le k \le \lfloor (B-1) / 2 \rfloor$
    \State $R \gets 0; \; c_3 \gets (\tfrac{1}{2})_{\lfloor (B-1) / 3 \rfloor} / (\lfloor (B-1) / 3 \rfloor)!$
    \For {$(m_3 \gets \lfloor (B-1) / 3 \rfloor; \; m_3 \ge 0; \; m_3 \gets m_3 - 1)$}
      \If {$m_3 \ne \lfloor (B-1) / 3 \rfloor$}
        \State $c_3 \gets c_3 \cdot (2m_3+2) / (2m_3+1)$
      \EndIf
      \State $s \gets 0; \; c_2 \gets c_3$
      \For {$(m_2 \gets 0; \; 2 m_2 + 3 m_2 < B; \; m_2 \gets m_2 + 1)$}
          \State $s \gets s + E_2^{m_2} \cdot (-1)^{m_2} c_2 / (4 m_2 + 6 m_3 + 1)$
          \State $c_2 \gets c_2 \cdot (2m_2 + 2m_3 + 1) / (2m_2 + 2)$
      \EndFor
      \State $R \gets (R \cdot E_3) + s$
    \EndFor
    \State \Return $R$
    \EndProcedure
    \State \Return $A^{-1/2} \left(\Call{RSum}{E_2,E_3,B} + [\pm \varepsilon] + [\pm \varepsilon]i\right)$ \Comment{Include prefactor and error bound}
  \end{algorithmic}
\end{algorithm}

\index{argument reduction}
\index{hypergeometric series}
\index{rectangular splitting}
\index{Arb library}

\begin{figure}
\begin{center}
\begin{tikzpicture}[scale=0.4]
\draw node[fill,circle,scale=0.4] at (0,0) {};
\draw node[fill,circle,scale=0.4] at (0,1) {};
\draw node[fill,circle,scale=0.4] at (0,2) {};
\draw node[fill,circle,scale=0.4] at (0,3) {};

\draw node[fill,circle,scale=0.4] at (1,0) {};
\draw node[fill,circle,scale=0.4] at (1,1) {};
\draw node[fill,circle,scale=0.4] at (1,2) {};

\draw node[fill,circle,scale=0.4] at (2,0) {};
\draw node[fill,circle,scale=0.4] at (2,1) {};

\draw node[fill,circle,scale=0.4] at (3,0) {};
\draw node[fill,circle,scale=0.4] at (3,1) {};

\draw node[fill,circle,scale=0.4] at (4,0) {};

\node at (6,0) {$m_2$};
\node at (-1,3.3) {$m_3$};


\draw[->,-latex] (0,0) -- (0,3.7);
\draw[->,-latex] (0,0) -- (5,0);
\end{tikzpicture}
\caption{Exponents $2 m_2 + 3 m_3 < B$ appearing for $B = 10$.}
\label{fig:lattice}
\end{center}
\end{figure}
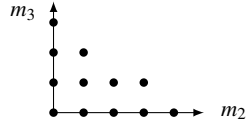

Algorithm~\ref{alg:hyprs} gives pseudocode for the method
implemented in Arb to compute $R_F$ using
combined argument reduction and series evaluation.
The subroutine RSum evaluates the series for $R_{-1/2}$
truncated to an arbitrary order $B$ using rectangular splitting
combined with recurrence relations for the coefficients
(one more optimization used in the implementation but omitted
from the pseudocode is to clear denominators so that all coefficients are small integers).
The exponents of $E_2^{m_2} E_3^{m_3}$ appearing in the series (Fig.~\ref{fig:lattice})
are the lattice points $m_2, m_3 \in \mathbb{Z}_{\ge 0}$ with $2m_2 + 3m_3 < B$:
we compute powers of $E_2$ and then use Horner's rule with respect to $E_3$.

We now consider the choice of $B$.
Since the series evaluation costs $O(B^2)$ operations and the argument
reduction costs $O(p / B)$ operations, the overall cost is
$O(N^2 + p/N)$ operations, which is minimized by setting $B \approx p^{1/3}$;
this gives us an $O(p^{1.667})$ bit complexity algorithm for
evaluating $R_F$.
With rectangular splitting for the series evaluation,
the optimal~$B$ should be closer to $B \approx p^{0.5}$ for moderate~$p$.
Timings in Arb
show that expanding to order $B = 2p^{0.4}$ for real variables and
$B = 2.5p^{0.4}$ for complex variables is nearly optimal
in the range $10 \le p \le 10^6$ bits.
Compared to the fixed order $B = 8$, this results in a measured speedup
of 1.5, 4.0, 11 and 31 times at a precision of 100, 1000, 10\,000 and 100\,000 decimal digits
respectively.

The algorithm for the $R_J$ series is essentially the same,
except that the summation
uses four nested loops instead of two to iterate
over the exponents with $2 m_2 + 3 m_3 + 4 m_4 + 5 m_5 < B$,
with corresponding nested recurrence relations to update
coefficients $c_2, c_3, c_4, c_5$ (see the Arb source code for details).
In this case,
rectangular splitting is used to split
the number of variables in half by precomputing a
twodimensional array of powers $E_2^{m_2} E_3^{m_3}$
and using the Horner scheme with respect to $E_4$ and $E_5$.
The speedup of combining rectangular splitting
with optimal argument reduction is smaller for $R_J$ than for $R_F$ but
still appreciable at very high precision.

\section{Arb implementation benchmarks}

\vspace{1mm}\noindent
Table~\ref{tab:timings} compares the performance of different
functions implemented in Arb.

\index{Arb library}
\index{elliptic integral}
\index{complete elliptic integral}
\index{incomplete elliptic integral}
\index{modular form}
\index{modular function}
\index{$j$-invariant}
\index{Dedekind eta function}
\index{Jacobi theta function}
\index{Weierstrass elliptic function}
\index{arithmetic-geometric mean}
\index{Carlson symmetric form}
\index{Legendre form}

The complete elliptic integrals of the first and second
kind are about as fast as the elementary functions at high precision
due to the $\OOtilde(p)$ AGM algorithm.\footnote{At precision up to about 1000 digits, the elementary functions in Arb are
significantly faster than the AGM due to using precomputed lookup tables and many low-level
optimizations~\cite{Johansson2015elementary}.}
The modular forms and functions which use $\OOtilde(p^{1.5})$ algorithms
with very low overhead are nearly as fast as complete elliptic integrals
in practice.

\begin{table}[t!]
\begin{centering}
\setlength{\tabcolsep}{0.5em}
\renewcommand{\arraystretch}{1.10}
\begin{tabular}{ l l | c c c c c }
& Function & $d=10$ & $d=10^2$ & $d=10^3$ & $d=10^4$ & $d=10^5$ \\ \hline
Elementary & $\exp(x)$ & $7.7 \cdot 10^{-7}$ & $2.9 \cdot 10^{-6}$ & $0.00011$ & $0.0062$ & $0.24$ \\
functions & $\log(x)$ & $8.1 \cdot 10^{-7}$ & $2.8 \cdot 10^{-6}$ & $0.00011$ & $0.0077$ & $0.27$ \\ \hline

Modular & $\eta(t)$ & $6.2 \cdot 10^{-6}$ & $1.99 \cdot 10^{-5}$ & $0.00037$ & $0.015$ & $0.69$ \\
forms \& & $j(t)$ & $6.3 \cdot 10^{-6}$ & $2.3 \cdot 10^{-5}$ & $0.00046$ & $0.022$ & $1.1$ \\
functions & $(\theta_i(0,t))_{i=1}^4$ & $7.6 \cdot 10^{-6}$ & $2.7 \cdot 10^{-5}$ & $0.00044$ & $0.022$ & $1.1$ \\ \hline

Elliptic & $(\theta_i(x,t))_{i=1}^4$ & $2.8 \cdot 10^{-5}$ & $8.1 \cdot 10^{-5}$ & $0.0016$ & $0.089$ & $5.4$ \\
and theta & $\wp(x,t)$ & $3.9 \cdot 10^{-5}$ & $0.00012$ & $0.0021$ & $0.11$ & $6.6$ \\
functions & $(\wp(x,t), \wp'(x,t))$ & $5.6 \cdot 10^{-5}$ & $0.00017$ & $0.0026$ & $0.13$ & $7.3$ \\
 & $\zeta(x,t)$ & $7.5 \cdot 10^{-5}$ & $0.00022$ & $0.0028$ & $0.14$ & $7.8$ \\
& $\sigma(x,t)$ & $7.6 \cdot 10^{-5}$ & $0.00022$ & $0.0030$ & $0.14$ & $8.1$ \\ \hline

Complete & $K(x)$ & $5.4 \cdot 10^{-6}$ & $2.0 \cdot 10^{-5}$ & $0.00018$ & $0.0068$ & $0.23$ \\
elliptic & $E(y)$ & $1.7 \cdot 10^{-5}$ & $6.1 \cdot 10^{-5}$ & $0.00072$ & $0.025$ & $0.71$ \\
integrals & $\Pi(x,y)$ & $7.0 \cdot 10^{-5}$ & $0.00046$ & $0.014$ & $3.6$ & $563$ \\ \hline

Incomplete & $\wp^{-1}(x,t)$ & $3.1 \cdot 10^{-5}$ & $0.00014$ & $0.0025$ & $0.20$ & $20$ \\
elliptic & $F(x,y)$ & $2.4 \cdot 10^{-5}$ & $0.00011$ & $0.0022$ & $0.19$ & $19$ \\
integrals & $E(x,y)$ & $5.6 \cdot 10^{-5}$ & $0.00030$ & $0.0070$ & $0.76$ & $97$ \\
& $\Pi(x,y,z)$ & $0.00017$ & $0.00098$ & $0.030$ & $5.6$ & $895$ \\
& $R_F(x,y,z)$ & $1.6 \cdot 10^{-5}$ & $9.5 \cdot 10^{-5}$ & $0.0020$ & $0.18$ & $18$ \\
& $R_G(x,y,z)$ & $4.7 \cdot 10^{-5}$ & $0.00027$ & $0.0067$ & $0.75$ & $95$ \\
& $R_D(x,y,z)$ & $2.1 \cdot 10^{-5}$ & $0.00016$ & $0.0046$ & $0.57$ & $78$ \\
& $R_J(x,y,z,t)$ & $3.4 \cdot 10^{-5}$ & $0.00031$ & $0.012$ & $2.6$ & $428$ \\
\end{tabular}
\caption{Time in seconds to evaluate the function (or tuple of function values simultaneously)
at $d$ decimal digits of precision ($p = \lceil d \log_2{10} \rceil$ bits) for $d$ between 10 and 100\,000.
The arguments are set to generic complex numbers $x = \sqrt{2} + \sqrt{3} i, y = \sqrt{3} + \sqrt{5} i, z = \sqrt{5} + \sqrt{7} i, t = \sqrt{7} + i/\sqrt{11}$.
}
\label{tab:timings}
\end{centering}
\end{table}

Elliptic functions and Jacobi theta functions, also implemented with $\OOtilde(p^{1.5})$
algorithms, are some 5-10 times slower than the special case of
theta constants or modular forms.
The incomplete elliptic integrals based on the $R_F$ function
implemented with $O(p^{1.667})$ complexity have similar
performance to the elliptic functions
at moderate precision with a slight divergence becoming visible only at several thousand digits.
Indeed, $\wp(x,t)$ and $\wp^{-1}(x,t)$ have virtually identical performance
although the algorithms are completely independent.

The incomplete elliptic integrals based on the $R_J$ function
stand out as being noticeably slower than the other functions,
as a result of the more complicated argument reduction
and high-dimensional series expansion.

\section{Other methods}
\vspace{1mm}\noindent
Many numerical techniques apart from those covered in this text
are useful in connection with elliptic functions and modular forms.
Without going into detail,
we sketch a few important ideas.

\subsection{Quadratically convergent methods and Newton iteration}

\label{sect:quadratic}

The algorithms described above for complete elliptic integrals
have quasioptimal $\OOtilde(p)$ bit complexity
owing to the quadratically convergent AGM iteration,
while the algorithms for all other functions
have $\OOtilde(p^{1.5})$ or worse bit complexity.
In fact, it is possible to compute general elliptic functions, modular forms
and incomplete elliptic integrals with $\OOtilde(p)$ bit complexity
using generalizations of the AGM iteration
together with Newton's method for inverse functions.
We have omitted these methods in the present work since they are more complicated,
especially for complex variables, and
not necessarily faster for~$p$ encountered in practice.

\index{Newton iteration}
\index{quadratic convergence}
\index{Landen transformation}
\index{elliptic logarithm}
\index{Jacobi theta function}
\index{modular form}
\index{elliptic integral}

The asymptotically fast computation of modular forms and modular functions
is discussed by Dupont~\cite{dupont2011fast}, and 
Labrande~\cite{labrande2017computing} has given algorithms
for general theta functions and elliptic functions.
An important special case is the inverse Weierstrass elliptic function
in the form of the elliptic logarithm, which
can be computed using a simple
AGM-type algorithm~\cite{cremona2013complex}.
For the Legendre incomplete elliptic integrals, algorithms based on
the quadratic Landen transformations are classical and have
been described in several other works; they have the disadvantage
of involving trigonometric functions, not having a straightforward
extension to complex variables, and in some regions
suffering from precision loss.

\subsection{Numerical integration}

Direct numerical integration is a viable way
to compute elliptic integrals.
Numerical integration is generally slower
than the more specialized algorithms already presented, but
with a robust general-purpose integration
algorithm, we can just plug in the formula for any
particular elliptic integral.
Specifying an explicit contour of integration
also provides full control over branch cuts.

\index{numerical integration}
\index{double exponential quadrature}

The double exponential or tanh-sinh quadrature method~\cite{takahasi1974double,bailey2011high}
is ideal for elliptic integrals since it is extremely simple and
converges rapidly even if the integrand
has algebraic singularities of unknown type at one or both endpoints.
The quadrature error in the double exponential method can be estimated quite reliably
using heuristics, and effective rigorous error bounds
are also known \cite{Molin2010a}.
Alternatively, Gauss-Jacobi quadrature can be used
for integrals with known algebraic singularities.
Recently, Molin and Neurohr have studied use of both
double exponential and Gauss-Jacobi quadrature
with rigorous error bounds for integration
of algebraic functions in the context of computing period
matrices for hyperelliptic curves~\cite{molin2017computing}.
Rigorous numerical integration code also exists in Arb~\cite{Johansson2018numerical},
but endpoint singularities require manual processing.

\index{trapezoidal rule}
\index{complete elliptic integral}

For integrals of smooth periodic functions, including integrals of
analytic functions on closed circular contours,
direct application of the trapezoidal rule
is often the best choice.
We conclude with the anecdote that Poisson already in the 1820s demonstrated the
use of the trapezoidal rule to approximate the
elliptic integral
$$\frac{1}{2\pi} \int_0^{2\pi} \sqrt{1 - 0.36 \sin^2(\theta)} d\theta$$
which is equal to $\tfrac{2}{\pi} E(0.36)$
in the Legendre notation.
Poisson derived an error bound for the $N$-point trapezoidal
approximation and showed that $N=16$ gives an error
less than $4.84 \cdot 10^{-6}$ for this integral (in fact, nine digits are correct).
Due to symmetry, just three nontrivial evaluations of the integrand
are required for this level of accuracy!
Trefethen and Weideman~\cite{trefethen2014exponentially} discuss
this example and provide a general error analysis for the
trapezoidal rule applied to periodic functions.




\vspace*{4mm}
\noindent
{\bf Acknowledgment.} The author thanks the organizers
of the KMPB Conference on \emph{Elliptic Integrals, Elliptic Functions and
Modular Forms in Quantum Field Theory}
for the invitation to present this work at DESY in October 2017
and for the opportunity to publish this extended
review in the post-conference proceedings.



\begin{thebibliography}{10}

\bibitem{bailey2011high}
D.~H. Bailey and J.~M. Borwein.
\newblock High-precision numerical integration: Progress and challenges.
\newblock {\em J. Symbolic Comput.}, 46(7):741--754, 2011.

\bibitem{bailey2015high}
D.~H. Bailey and J.~M. Borwein.
\newblock High-precision arithmetic in mathematical physics.
\newblock {\em Mathematics}, 3(2):337--367, 2015.

\bibitem{brent2010modern}
R.~P. Brent and P.~Zimmermann.
\newblock {\em Modern Computer Arithmetic}.
\newblock Cambridge University Press, 2010.

\bibitem{carlson1995numerical}
B.~C. Carlson.
\newblock Numerical computation of real or complex elliptic integrals.
\newblock {\em Numerical Algorithms}, 10(1):13--26, 1995.

\bibitem{cohen2013course}
H.~Cohen.
\newblock {\em A course in computational algebraic number theory}.
\newblock Springer Science \& Business Media, 2013.

\bibitem{cox2000arithmetic}
D.~A. Cox.
\newblock The arithmetic-geometric mean of {G}auss.
\newblock In {\em Pi: A source book}, pages 481--536. Springer, 2000.

\bibitem{cremona2013complex}
J.~E. Cremona and T.~Thongjunthug.
\newblock The complex {AGM}, periods of elliptic curves over {C} and complex
  elliptic logarithms.
\newblock {\em Journal of Number Theory}, 133(8):2813--2841, 2013.

\bibitem{dupont2006moyenne}
R.~Dupont.
\newblock {\em Moyenne arithm{\'e}tico-g{\'e}om{\'e}trique, suites de
  {B}orchardt et applications}.
\newblock PhD thesis, {\'E}cole polytechnique, Palaiseau, 2006.

\bibitem{dupont2011fast}
R.~Dupont.
\newblock Fast evaluation of modular functions using {N}ewton iterations and
  the {AGM}.
\newblock {\em Mathematics of Computation}, 80(275):1823--1847, 2011.

\bibitem{enge2009complexity}
A.~Enge.
\newblock The complexity of class polynomial computation via floating point
  approximations.
\newblock {\em Mathematics of Computation}, 78(266):1089--1107, 2009.

\bibitem{enge2018short}
A.~Enge, W.~Hart, and F.~Johansson.
\newblock Short addition sequences for theta functions.
\newblock {\em Journal of Integer Sequences}, 21(2):3, 2018.

\bibitem{Fieker2017nemo}
C.~Fieker, W.~Hart, T.~Hofmann, and F.~Johansson.
\newblock {Nemo/Hecke}: computer algebra and number theory packages for the
  {J}ulia programming language.
\newblock In {\em Proceedings of the 42nd International Symposium on Symbolic
  and Algebraic Computation}, ISSAC '17, pages 1--1, Kaiserslautern, Germany,
  2017. ACM.

\bibitem{izzo2016astrodynamics}
D.~Izzo and F.~Biscani.
\newblock On the astrodynamics applications of {W}eierstrass elliptic and
  related functions.
\newblock \url{https://arxiv.org/abs/1601.04963}, 2016.

\bibitem{Johansson2015elementary}
F.~Johansson.
\newblock Efficient implementation of elementary functions in the
  medium-precision range.
\newblock In {\em 22nd IEEE Symposium on Computer Arithmetic}, ARITH22, pages
  83--89, 2015.

\bibitem{Johansson2016hypergeometric}
F.~Johansson.
\newblock Computing hypergeometric functions rigorously, 2016.
\newblock \url{http://arxiv.org/abs/1606.06977}.

\bibitem{Johansson2017arb}
F.~Johansson.
\newblock Arb: efficient arbitrary-precision midpoint-radius interval
  arithmetic.
\newblock {\em IEEE Transactions on Computers}, 66:1281--1292, 2017.

\bibitem{Johansson2018numerical}
F.~Johansson.
\newblock Numerical integration in arbitrary-precision ball arithmetic, 2018.
\newblock \url{https://arxiv.org/abs/1802.07942}.

\bibitem{labrande2017computing}
H.~Labrande.
\newblock Computing {J}acobi's theta in quasi-linear time.
\newblock {\em Mathematics of Computation}, 2017.

\bibitem{Molin2010a}
P.~Molin.
\newblock {\em Numerical Integration and {L}-Functions Computations}.
\newblock Theses, Universit{\'e} Sciences et Technologies - Bordeaux I, October
  2010.

\bibitem{molin2017computing}
P.~Molin and C.~Neurohr.
\newblock Computing period matrices and the {A}bel-{J}acobi map of
  superelliptic curves.
\newblock {\em arXiv preprint arXiv:1707.07249}, 2017.

\bibitem{muller2006elementary}
J.~M. Muller.
\newblock {\em Elementary Functions}.
\newblock Springer, 2006.

\bibitem{nogneng2018evaluation}
D.~Nogneng and E.~Schost.
\newblock On the evaluation of some sparse polynomials.
\newblock {\em Math. Comp}, 87:893--904, 2018.

\bibitem{Olver2010}
F.~W.~J. Olver, D.~W. Lozier, R.~F. Boisvert, and C.~W. Clark.
\newblock {\em {NIST Handbook of Mathematical Functions}}.
\newblock Cambridge University Press, New York, 2010.

\bibitem{PatersonStockmeyer1973}
M.~S. Paterson and L.~J. Stockmeyer.
\newblock On the number of nonscalar multiplications necessary to evaluate
  polynomials.
\newblock {\em SIAM Journal on Computing}, 2(1), March 1973.

\bibitem{Rademacher1973}
H.~Rademacher.
\newblock {\em Topics in analytic number theory}.
\newblock Springer, 1973.

\bibitem{Smith1989}
D.~M. Smith.
\newblock Efficient multiple-precision evaluation of elementary functions.
\newblock {\em Mathematics of Computation}, 52:131--134, 1989.

\bibitem{takahasi1974double}
H.~Takahasi and M.~Mori.
\newblock Double exponential formulas for numerical integration.
\newblock {\em Publ. Res. Inst. Math. Sci.}, 9(3):721--741, 1974.

\bibitem{PARI2}
{The PARI~Group}, Univ. Bordeaux.
\newblock {\em {PARI/GP version 2.9.4}}, 2017.

\bibitem{sagemath}
{The Sage Developers}.
\newblock {\em {S}ageMath, the {S}age {M}athematics {S}oftware {S}ystem
  ({V}ersion 8.2.0)}, 2018.
\newblock {\tt http://www.sagemath.org}.

\bibitem{trefethen2014exponentially}
L.~N. Trefethen and J.~A.~C. Weideman.
\newblock The exponentially convergent trapezoidal rule.
\newblock {\em SIAM Review}, 56(3):385--458, 2014.

\bibitem{wolfellip}
{W}olfram {R}esearch.
\newblock The {W}olfram {F}unctions {S}ite - {E}lliptic {I}ntegrals.
\newblock \url{http://functions.wolfram.com/EllipticIntegrals/}, 2016.

\end{thebibliography}
\end{document}